\documentclass[reqno]{amsart}
\usepackage{setspace,tikz,xcolor,mathrsfs,listings,multicol}
\usepackage{rotating}
\usepackage[vcentermath]{youngtab}
\usepackage[margin=1in,includefoot,footskip=30pt]{geometry}
\usepackage{enumerate}
\usepackage{booktabs}
\usepackage[all,cmtip]{xy}
\usetikzlibrary{arrows,matrix}
\tikzset{tab/.style={matrix of math nodes,column sep=-.35, row sep=-.35,text height=7pt,text width=7pt,align=center,inner sep=2,font=\footnotesize}}

\usepackage[colorlinks=true, pdfstartview=FitV, linkcolor=blue, citecolor=blue, urlcolor=blue]{hyperref}

\newcommand{\g}{\mathfrak{g}}
\newcommand{\clfw}{\overline{\Lambda}} 
\newcommand{\inner}[2]{\left\langle #1, #2 \right\rangle}
\newcommand{\iso}{\cong}

\newcommand{\qbinom}[3]{\genfrac{[}{]}{0pt}{}{#1}{#2}_{#3}}

\newcommand{\Absval}[1]{\left\lVert #1 \right\rVert}

\newcommand{\te}{\widetilde{e}}
\newcommand{\tf}{\widetilde{f}}

\DeclareMathOperator{\wt}{wt} 

\newcommand{\mcA}{\mathcal{A}}

\newcommand{\ZZ}{\mathbb{Z}}
\newcommand{\QQ}{\mathbb{Q}}

\definecolor{darkred}{rgb}{0.7,0,0} 
\newcommand{\defn}[1]{{\color{darkred}\emph{#1}}} 

\definecolor{UQgold}{RGB}{196, 158, 54} 
\definecolor{UQpurple}{RGB}{73, 7, 94} 
\definecolor{UMNgold}{RGB}{255,200,46} 
\definecolor{UMNmaroon}{RGB}{106,0,50} 

\usepackage{listings}
\lstdefinelanguage{Sage}[]{Python}
{morekeywords={False,sage,True},sensitive=true}
\definecolor{dblackcolor}{rgb}{0.0,0.0,0.0}
\definecolor{dbluecolor}{rgb}{0.01,0.02,0.7}
\definecolor{dgreencolor}{rgb}{0.2,0.4,0.0}
\definecolor{dgraycolor}{rgb}{0.30,0.3,0.30}

\usepackage{xparse}

\makeatletter
\protected\def\specialmergetwolists{%
  \begingroup
  \@ifstar{\def\cnta{1}\@specialmergetwolists}
    {\def\cnta{0}\@specialmergetwolists}%
}
\def\@specialmergetwolists#1#2#3#4{%
  \def\tempa##1##2{%
    \edef##2{%
      \ifnum\cnta=\@ne\else\expandafter\@firstoftwo\fi
      \unexpanded\expandafter{##1}%
    }%
  }%
  \tempa{#2}\tempb\tempa{#3}\tempa
  \def\cnta{0}\def#4{}%
  \foreach \x in \tempb{%
    \xdef\cnta{\the\numexpr\cnta+1}%
    \gdef\cntb{0}%
    \foreach \y in \tempa{%
      \xdef\cntb{\the\numexpr\cntb+1}%
      \ifnum\cntb=\cnta\relax
        \xdef#4{#4\ifx#4\empty\else,\fi\x#1\y}%
        \breakforeach
      \fi
    }%
  }%
  \endgroup
}
\makeatother

\usepackage{xparse}
\DeclareDocumentCommand\rpp{ m m g }{
	\foreach \x [count=\s from 1] in {#1}{
	        {\ifnum\s=1
	                \draw (0,-\s)--(\x,-\s);
	                \fi}
	   \draw (0,-\s-1) to (\x,-\s-1);
	   \foreach \y in {0, ..., \x} {\draw (\y,-\s)--(\y,-\s-1);}
	}
	\specialmergetwolists{/}{#1}{#2}\ziplist
	\foreach \x/\y [count=\yi from 1] in \ziplist{
	    \node[anchor=west,font=\small] at (\x,-\yi - .5) {$\y$};
	}
	\IfValueT {#3}
	{\foreach \z [count=\zi from 1] in {#3} {\node[anchor=east,font=\small] at (0,-\zi - .5) {$\z$};}}
	{}
}

\theoremstyle{plain}
\newtheorem{theorem}{Theorem}[section]
\newtheorem{lemma}[theorem]{Lemma}

\newtheorem{proposition}[theorem]{Proposition}

\theoremstyle{definition}

\numberwithin{equation}{section}

\usepackage[colorinlistoftodos]{todonotes}

\begin{document}
\title[Existence KR crystals when multiplicity free]{Existence of Kirillov--Reshetikhin crystals for multiplicity free nodes}

\author[R.~Biswal]{Rekha Biswal}
\address[R. Biswal]{Max-Planck-Institut f\"ur Mathematik, Vivatsgasse 7, 53111 Bonn, Germany}
\email{rekha@mpim-bonn.mpg.de}
\urladdr{https://rekhabiswal.github.io/}

\author[T.~Scrimshaw]{Travis Scrimshaw}
\address[T.~Scrimshaw]{School of Mathematics and Physics, The University of Queensland, St.\ Lucia, QLD 4072, Australia}
\email{tcscrims@gmail.com}
\urladdr{https://people.smp.uq.edu.au/TravisScrimshaw/}


\thanks{T.S.~was partially supported by the Australian Research Council DP170102648.}

\begin{abstract}
We show that the Kirillov--Reshetikhin crystal $B^{r,s}$ exists when $r$ is a node such that the Kirillov--Reshetikhin module $W^{r,s}$ has a multiplicity free classical decomposition.
\end{abstract}

\maketitle

\section{Introduction}
\label{sec:introduction}

Kirillov--Reshetikhin (KR) modules are an class of finite-dimensional representation of an affine quantum group $U_q'(\g)$ without the degree operator that is classified by their Drinfel'd polynomials that have received significant attention.
We denote a KR module by $W^{r,s}$, where $r$ is a node of the classical (\textit{i.e.} underlying finite type) Dynkin diagram and $s \in \ZZ_{>0}$. One construction of a KR module $W^{r,s}$ is by computing the minimal affinization of the highest weight $U_q(\g_0)$-module $V(s\clfw_r)$~\cite{Chari95,CP95II,CP96,CP96II}, where $\g_0$ is the classical Lie algebra. Another method is by using the fusion construction of~\cite{KKMMNN92} from the image under an $R$-matrix of an $s$-fold tensor product of the fundamental module $W^{r,1}$ (see, \textit{e.g.},~\cite{Kashiwara02}). KR modules are also known to have special properties. The classical decomposition, the branching rule of $W^{r,s}$ to a $U_q(\g_0)$-module, is given by a fermionic formula~\cite{dFK08,Hernandez10}, which leads to the (virtual) Kleber algorithm~\cite{Kleber98,OSS03II}.
The characters (resp.\ $q$-characters) of KR modules also satisfy the $Q$-system (resp.\ $T$-system) relations~\cite{Hernandez10,Nakajima03II}.
Furthermore, the graded characters of (Demazure submodules of) a tensor product of fundamental modules are (nonsymmetric) Macdonald polynomials at $t=0$~\cite{LNSSS14,LNSSS14II} (\cite{LNSSS15}).

One important (conjectural) property~\cite{HKOTY99,HKOTT02} is that the KR module $W^{r,s}$ admits a crystal base~\cite{K90,K91}, which is known as a Kirillov--Reshetikhin (KR) crystal and denoted by $B^{r,s}$.
Kashiwara showed that all fundamental modules $W^{r,1}$ have crystal bases~\cite{Kashiwara02}.
It was shown that $B^{r,s}$ exists in all nonexceptional types in~\cite{Okado07,OS08} and in types $G_2^{(1)}$ and $D_4^{(3)}$ in~\cite{KMOY07,Naoi17,Yamane98}. For all affine types, the existence of $B^{r,s}$ has been proven when $r$ is adjacent to $0$ or in the orbit of $0$ under a Dynkin diagram automorphism (equivalently, $W^{r,s}$ is irreducible as $U_q(\g)$-module)~\cite{KKMMNN92}.

Our main result is that the KR module $W^{r,s}$ has a crystal base whenever its classical decomposition is multiplicity free in all affine types. We do this by showing the existence of $B^{r,s}$ in the cases not covered by~\cite{KKMMNN92,Okado07,OS08}. More explicitly, we show this for $r = 3,5$ in type $E_6^{(1)}$, for $r=2,6$ in type $E_7^{(1)}$, for $r=7$ in type $E_8^{(1)}$, and for $r=4$ in types $F_4^{(1)}$ and $E_6^{(2)}$, where we label the Dynkin diagrams following~\cite{Bourbaki02} (see also Figure~\ref{fig:exceptional_types} for the labeling).
Using the techniques developed in~\cite{KKMMNN92}, our proof shows the existence of a crystal pseudobase $(L, B)$ by using the fusion construction of $W^{r,s}$ and is similar to~\cite{Okado07,OS08} by calculating the prepolarization for certain vectors. From there, we can construct the associated crystal by $B / \{\pm1\}$.

Let us describe some possible applications of our results.
The $X = M$ conjecture~\cite{HKOTY99,HKOTT02} arises from mathematical physics relating vertex models and the Bethe ansatz of Heisenberg spin chains, and the $X$ side requires the existence of KR crystals.
A uniform model for $B^{r,1}$ was given using quantum and projected level-zero LS paths~\cite{LNSSS14,LNSSS16,LNSSS14II,NS06II,NS08II,NS08}.
Since the KR crystal $B^{r,s}$ exists, we have a partial (conjectural) combinatorial description from~\cite{LS18} using $(B^{r,1})^{\otimes s}$, partially mimicking the fusion construction.

After completion of this paper, we learned that Naoi independently proved all cases in type $E_6^{(1)}$ using similar techniques~\cite{Naoi19}.

This paper is organized as follows.
In Section~\ref{sec:background}, we give the necessary background.
In Section~\ref{sec:existence}, we show our main result: that the KR modules $W^{r,s}$ has a crystal pseudobase whenever $W^{r,s}$ has a multiplicity free classical decomposition.

\subsection*{Acknowledgements}

The authors would like to thank Katsuyuki Naoi and Masato Okado for useful discussions.
RB would like to thank Tokyo University of Agriculture and Technology for its hospitality during her visit in Ju
This work benefited from computations using {\sc SageMath}~\cite{sage,combinat}.

\section{Background}
\label{sec:background}

In this section, we provide the necessary background.

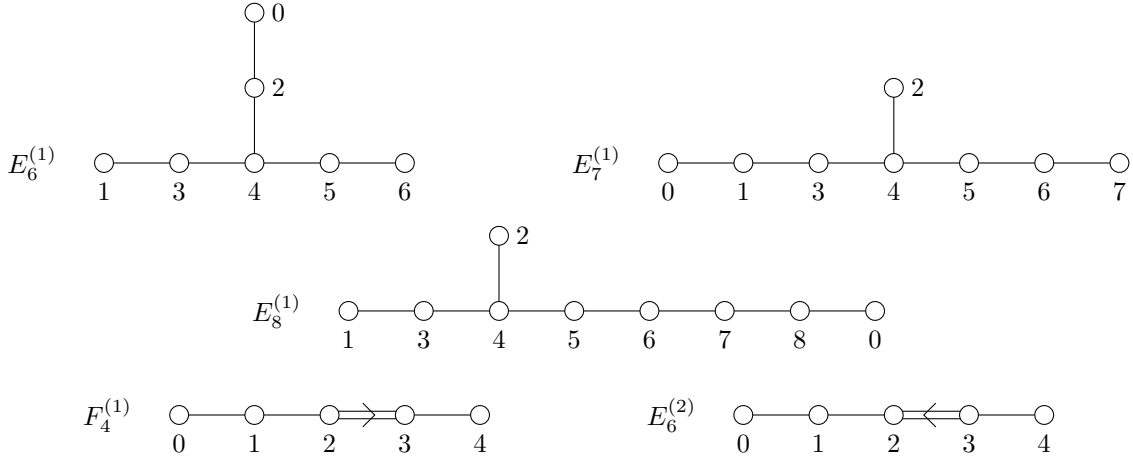
\begin{figure}[t]
\begin{center}
\begin{tikzpicture}[scale=0.5, baseline=-20]
\draw (-1,0) node[anchor=east] {$E_6^{(1)}$};
\draw (0 cm,0) -- (8 cm,0);
\draw (4 cm, 0 cm) -- +(0,2 cm);
\draw (4 cm, 2 cm) -- +(0,2 cm);
\draw[fill=white] (4 cm, 4 cm) circle (.25cm) node[right=3pt]{$0$};
\draw[fill=white] (0 cm, 0 cm) circle (.25cm) node[below=4pt]{$1$};
\draw[fill=white] (2 cm, 0 cm) circle (.25cm) node[below=4pt]{$3$};
\draw[fill=white] (4 cm, 0 cm) circle (.25cm) node[below=4pt]{$4$};
\draw[fill=white] (6 cm, 0 cm) circle (.25cm) node[below=4pt]{$5$};
\draw[fill=white] (8 cm, 0 cm) circle (.25cm) node[below=4pt]{$6$};
\draw[fill=white] (4 cm, 2 cm) circle (.25cm) node[right=3pt]{$2$};
\begin{scope}[xshift=15cm]
\draw (-1,0) node[anchor=east] {$E_7^{(1)}$};
\draw (0 cm,0) -- (12 cm,0);
\draw (6 cm, 0 cm) -- +(0,2 cm);
\draw[fill=white] (0 cm, 0 cm) circle (.25cm) node[below=4pt]{$0$};
\draw[fill=white] (2 cm, 0 cm) circle (.25cm) node[below=4pt]{$1$};
\draw[fill=white] (4 cm, 0 cm) circle (.25cm) node[below=4pt]{$3$};
\draw[fill=white] (6 cm, 0 cm) circle (.25cm) node[below=4pt]{$4$};
\draw[fill=white] (8 cm, 0 cm) circle (.25cm) node[below=4pt]{$5$};
\draw[fill=white] (10 cm, 0 cm) circle (.25cm) node[below=4pt]{$6$};
\draw[fill=white] (12 cm, 0 cm) circle (.25cm) node[below=4pt]{$7$};
\draw[fill=white] (6 cm, 2 cm) circle (.25cm) node[right=3pt]{$2$};
\end{scope}
\end{tikzpicture}

\begin{tikzpicture}[scale=0.5]
\draw (-1,0) node[anchor=east] {$E_8^{(1)}$};
\draw (0 cm,0) -- (12 cm,0);
\draw (4 cm, 0 cm) -- +(0,2 cm);
\draw (12 cm,0) -- +(2 cm,0);
\draw[fill=white] (14 cm, 0 cm) circle (.25cm) node[below=4pt]{$0$};
\draw[fill=white] (0 cm, 0 cm) circle (.25cm) node[below=4pt]{$1$};
\draw[fill=white] (2 cm, 0 cm) circle (.25cm) node[below=4pt]{$3$};
\draw[fill=white] (4 cm, 0 cm) circle (.25cm) node[below=4pt]{$4$};
\draw[fill=white] (6 cm, 0 cm) circle (.25cm) node[below=4pt]{$5$};
\draw[fill=white] (8 cm, 0 cm) circle (.25cm) node[below=4pt]{$6$};
\draw[fill=white] (10 cm, 0 cm) circle (.25cm) node[below=4pt]{$7$};
\draw[fill=white] (12 cm, 0 cm) circle (.25cm) node[below=4pt]{$8$};
\draw[fill=white] (4 cm, 2 cm) circle (.25cm) node[right=3pt]{$2$};
\end{tikzpicture}

\begin{tikzpicture}[scale=0.5,baseline=10]
\draw (-1,0) node[anchor=east] {$F_4^{(1)}$};
\draw (0 cm,0) -- (2 cm,0);
{
\pgftransformxshift{2 cm}
\draw (0 cm,0) -- (2 cm,0);
\draw (2 cm, 0.1 cm) -- +(2 cm,0);
\draw (2 cm, -0.1 cm) -- +(2 cm,0);
\draw (4.0 cm,0) -- +(2 cm,0);
\draw[shift={(3.2, 0)}, rotate=0] (135 : 0.45cm) -- (0,0) -- (-135 : 0.45cm);
\draw[fill=white] (0 cm, 0 cm) circle (.25cm) node[below=4pt]{$1$};
\draw[fill=white] (2 cm, 0 cm) circle (.25cm) node[below=4pt]{$2$};
\draw[fill=white] (4 cm, 0 cm) circle (.25cm) node[below=4pt]{$3$};
\draw[fill=white] (6 cm, 0 cm) circle (.25cm) node[below=4pt]{$4$};
}
\draw[fill=white] (0 cm, 0 cm) circle (.25cm) node[below=4pt]{$0$};

\begin{scope}[xshift=15cm]
\draw (-1,0) node[anchor=east] {$E_6^{(2)}$};
\draw (0 cm,0) -- (2 cm,0);
{
\pgftransformxshift{2 cm}
\draw (0 cm,0) -- (2 cm,0);
\draw (2 cm, 0.1 cm) -- +(2 cm,0);
\draw (2 cm, -0.1 cm) -- +(2 cm,0);
\draw (4.0 cm,0) -- +(2 cm,0);
\draw[shift={(2.8, 0)}, rotate=180] (135 : 0.45cm) -- (0,0) -- (-135 : 0.45cm);
\draw[fill=white] (0 cm, 0 cm) circle (.25cm) node[below=4pt]{$1$};
\draw[fill=white] (2 cm, 0 cm) circle (.25cm) node[below=4pt]{$2$};
\draw[fill=white] (4 cm, 0 cm) circle (.25cm) node[below=4pt]{$3$};
\draw[fill=white] (6 cm, 0 cm) circle (.25cm) node[below=4pt]{$4$};
}
\draw[fill=white] (0 cm, 0 cm) circle (.25cm) node[below=4pt]{$0$};
\end{scope}
\end{tikzpicture}
\end{center}
\caption{Dynkin diagrams for affine type $E_{6,7,8}^{(1)}$, $F_4^{(1)}$, and $E_6^{(2)}$.}
\label{fig:exceptional_types}
\end{figure}

Let $\g$ be an affine Kac--Moody Lie algebra with index set $I$, Cartan matrix $A = (A_{ij})_{i,j \in I}$, simple roots $(\alpha_i)_{i \in I}$, simple coroots $(h_i)_{i \in I}$, fundamental weights $(\Lambda_i)_{i \in I}$, weight lattice $P$, dominant weights $P^+$, coweight lattice $P^{\vee}$, and canonical pairing $\langle\ ,\ \rangle \colon P^{\vee} \times P \to \ZZ$ given by $\inner{h_i}{\alpha_j} = A_{ij}$.
We note that we follow the labeling given in~\cite{Bourbaki02} (see Figure~\ref{fig:exceptional_types} for the exceptional types and their labellings).
Let $\g_0$ denote the canonical simple Lie algebra given by the index set $I_0 = I \setminus \{0\}$.
Let $\overline{\lambda}$ denote the natural projection of $\lambda \in P$ onto the weight lattice $P_0$ of $\g_0$, so $\{\clfw_r\}_{r \in I_0}$ are the fundamental weights of $\g_0$.
Let $\varpi_r = \Lambda_r - \inner{c}{\Lambda_r} \Lambda_0$, where $c$ is the canonical central element of $\g$, denote the level-zero fundamental weights.
Let $q$ be an indeterminate, and we denote
\[
[m]_q = \dfrac{q^m - q^{-m}}{q-q^{-1}},
\qquad\qquad
[k]_q! = [k]_q [k-1]_q \dotsm [1]_q,
\qquad\qquad
\qbinom{m}{k}{q} = \dfrac{[m]_q[m-1]_q \dotsm [m-k+1]_q}{[k]_q!},
\]
for $m \in \ZZ$ and $k \in \ZZ_{\geq 0}$.
Let $q_i = q^{s_i}$ and $K_i = q^{s_i h_i}$, where $(s_1, \dotsc, s_n)$ is the diagonal symmetrizing matrix of $A$.

\subsection{Quantum groups}

Let $U_q'(\g) = U_q([\g, \g])$ denote the quantum group of the derived subalgebra of $\g$.
More specifically, the quantum group $U_q'(\g)$ is the associative $\QQ(q)$-algebra generated by $e_i, f_i, q^h$, where $i \in I$ and $h \in P^{\vee}$, that satisfies the relations
\begin{gather*}
q^0 = 1, \qquad q^h q^{h'} = q^{h+h'}, \qquad \text{for } h, h' \in P^{\vee},
\\ q^h e_i q^{-h} = q^{\inner{h}{\alpha_i}} e_i, \qquad q^h f_i q^{-h} = q^{-\inner{h}{\alpha_i}} f_i, \qquad \text{for } h \in P^{\vee}, i \in I,
\\ e_i f_j - f_j e_i = \delta_{ij} \dfrac{K_i - K_i^{-1}}{q_i - q_i^{-1}} \qquad \text {for } i,j \in I,
\end{gather*}
and the \defn{(quantum) Serre relations}
\[
\sum_{k=0}^{1-A_{ij}} (-1)^{k} e_i^{(k)} e_j e_i^{(1-A_{ij}-k)} = 0,
\qquad\qquad
\sum_{k=0}^{1-A_{ij}} (-1)^{k} f_i^{(k)} f_j f_i^{(1-A_{ij}-k)} = 0,
\]
where $e_i^{(k)} = e_i^k / [k]_{q_i}!$ and $f_i^{(k)} = f_i^k / [k]_{q_i}!$, for all $i,j \in I$ such that $i \neq j$.
We recall that $U_q'(\g)$ is a Hopf algebra; in particular, there exists a coproduct so we can take tensor products of $U_q'(\g)$-modules.

Denote the weight lattice of $U_q'(\g)$ by $P' = P / \ZZ \delta$, where $\delta$ is the null root of $\g$. Therefore, there is a linear dependence relation on the simple roots in $P'$. As we will not be considering $U_q(\g)$-modules in this paper, we will abuse notation and denote the $U_q'(\g)$-weight lattice by $P$.
For a $U_q'(\g)$-module $M$ and $\lambda \in P$, we denote the $\lambda$ weight space by
\[
M_{\lambda} = \{ v \in M \mid q^h v = q^{\inner{h}{\lambda}} v \text{ for all } h \in P^{\vee} \}.
\]
If $v \in M_{\lambda} \setminus \{0\}$, then we say $\wt(v) = \lambda$.

For $\lambda \in P_0^+$, we denote the highest weight $U_q(\g_0)$-module by $V(\lambda)$.

\subsection{Crystal (pseudo)bases and polarizations}

Let $\mcA$ denote the subring of $\QQ(q)$ of rational functions without poles at $0$.
A \defn{crystal base} of an integrable $U_q'(\g)$-module $M$ is a pair $(L, B)$, where $L$ is a free $\mcA$-module and $B$ is a basis of the $\QQ$-vector space $L / qL$, such that
\begin{enumerate}
\item $M \iso \QQ(q) \otimes_{\mcA} L$,
\item $L \iso \bigoplus_{\lambda \in P} L_{\lambda}$ with $L_{\lambda} = L \cap M_{\lambda}$,
\item $\te_i L \subseteq L$ and $\tf_i L \subseteq L$ for all $i \in I$,
\item $B = \bigsqcup_{\lambda \in P} B_{\lambda}$ with $B_{\lambda} = B \cap (L_{\lambda} / q L_{\lambda})$,
\item $\te_i B \subseteq B \sqcup \{0\}$ and $\tf_i B \subseteq B \sqcup \{0\}$,
\item $\tf_i b = b'$ if and only if $\te_i b' = b$ for all $b,b' \in B$ and $i \in I$.
\end{enumerate}
We say $(L, B)$ is a \defn{crystal pseudobase} of $M$ if it satisfies the conditions above for $B = B' \sqcup (-B')$, where $B'$ is a basis of $L / qL$.

For a $U_q'(\g)$-module $M$, a \defn{prepolarization} is a symmetric bilinear form $(\ ,\ ) \colon M \times M \to \QQ(q)$ that satisfies
\begin{equation}
\label{eq:admissible_pairing}
(q^h v, w) = (v, q^h w),
\qquad\qquad
(e_i v, w) = (v, q_i^{-1} K_i^{-1} f_i w),
\qquad\qquad
(f_i v, w) = (v, q_i^{-1} K_i e_i w),
\end{equation}
for all $i \in I$, and $v,w \in M$.\footnote{For $U_q(\g)$-modules $M, N$, a pairing $(\ ,\ ) \colon M \times N \to \QQ(q)$ that satisfies~\eqref{eq:admissible_pairing} is often called admissible.} Denote $\Absval{v}^2 = (v, v)$. If a prepolarization is positive definite with respect to the total order on $\QQ(q)$
\[
f > g \text{ if and only if } f - g \in \bigsqcup_{n \in \ZZ} \{ q^n(d + q\mcA) \mid d \in \QQ_{>0} \}
\]
(with $f \geq g$ defined as $f = g$ or $f > g$) then it is called a \defn{polarization}.

\subsection{Kirillov--Reshetikhin modules and the fusion construction}

Consider the subalgebras of $\QQ(q)$
\[
\mcA_{\ZZ} = \{ f(q) / g(q) \mid f(q), g(q) \in \ZZ[q], g(0) = 1 \},
\qquad\qquad
K_{\ZZ} = \mcA_{\ZZ}[q^{-1}].
\]
Let $U_q'(\g)_{K_{\ZZ}}$ denote the $K_{\ZZ}$-subalgebra of $U_q'(\g)$ generated by $e_i, f_i, q^h$ for all $i \in I$ and $h \in P^{\vee}$. The following is a combination of~\cite[Prop.~2.6.1]{KKMMNN92} and~\cite[Prop.~2.6.2]{KKMMNN92}.

\begin{proposition}
\label{prop:existence}
Let $M$ be a finite-dimensional integrable $U_q'(\g)$-module. Suppose $M$ has a prepolarization $(\ ,\ )$ and a $U_q'(\g)_{K_{\ZZ}}$-submodule $M_{K_{\ZZ}}$ such that $(M_{K_{\ZZ}}, M_{K_{\ZZ}}) \subseteq K_{\ZZ}$. Assume $M \iso \bigoplus_{k=1}^m V(\overline{\lambda}_k)$ as $U_q(\g_0)$-modules, with $\overline{\lambda}_k \in P_0^+$ for all $k$, such that there exists $u_k \in (M_{K_{\ZZ}})_{\lambda_k}$ such that $(u_k, u_{\ell}) \in \delta_{k\ell} + q\mcA$ and $\Absval{e_i u_k}^2 \in q_i^{-2\inner{h_i}{\lambda_k}-2} q \mcA$. Then $(\ ,\ )$ is a polarization and for
\[
L = \{v \in M \mid \Absval{v}^2 \in \mcA \}, \qquad B = \{ b \in (M_{K_{\ZZ}} \cap L) / (M_{K_{\ZZ}} \cap qL) \mid (b,b)_0 = 1 \},
\]
where $(\ ,\ )_0 \colon L/qL \to \QQ$ is the bilinear form induced by $(\ ,\ )$, the pair $(L, B)$ is a crystal pseudobase of $M$.
\end{proposition}

For an indeterminate $z$, let $M_z$ denote the $U_q'(\g)$-module $\QQ(q)[z,z^{-1}] \otimes M$, where $e_i$ and $f_i$ act by $z^{\delta_{0i}} \otimes e_i$ and $z^{-\delta_{0i}} \otimes f_i$ called the \defn{affinization module} of $M$.
For $a \in \QQ(q)$, define the \defn{evaluation module} $M_a = M_z / (z - a) M_z$.
For $v \in M$, let $v_a$ denote the corresponding element in $M_a$ (\textit{i.e.}, the projection of $1 \otimes v$).
Let $W(\varpi_r)$ denote the \defn{fundamental module} from~\cite{Kashiwara02}.

\begin{proposition}[{\cite[Prop.~9.3]{Kashiwara02}}]
Consider nonzero $a,b \in \QQ(q)$ such that $a/b \in \mcA$. Then for any $r \in I_0$, there exists a unique nonzero $U_q'(\g)$-module homomorphism
\[
R_{a,b} \colon W(\varpi_r)_a \otimes W(\varpi_r)_b \to W(\varpi_r)_b \otimes W(\varpi_r)_a
\]
that satisfies $R_{a,b}( u_a \otimes u_b ) = u_b \otimes u_a$ for some nonzero $u \in W(\varpi_r)_{\varpi_r}$. The map $R_{a,b}$ is called the \defn{(normalized) $R$-matrix} and satisfies the Yang--Baxter equation.
\end{proposition}

Denote
\[
W(\varpi_r; a_1, a_2, \dotsc, a_m) = W(\varpi_r)_{a_1} \otimes W(\varpi_r)_{a_2} \otimes \cdots \otimes W(\varpi_r)_{a_m}.
\]
Let $\kappa = s_i$ if $\g$ is of untwisted affine type and $\kappa = 1$  if $\g$ is of twisted affine type.
Since the $R$-matrix satisfies the Yang--Baxter equation, we can define the map
\[
R_s \colon W(\varpi_r; q^{\kappa(s-1)}, q^{\kappa(s-3)}, \dotsc, q^{\kappa(1-s)}) \to W(\varpi_r; q^{\kappa(1-s)}, \dotsc, q^{\kappa(s-3)}, q^{\kappa(s-1)})
\]
by applying the $R$-matrix on every pair of factors according to the long element of the symmetric group on $s$ letters $(q^{\kappa(s-1)}, q^{\kappa(s-3)}, \dotsc, q^{\kappa(1-s)})$.
Let $W^{r,s}$ denote the image of $R_s$, which is a simple $U_q'(\g)$-module~\cite{Kashiwara02}, and we call $W^{r,s}$ a \defn{Kirillov--Reshetikhin (KR) module}. From~\cite{CP95,CP98}, the module $W^{r,s}$ satisfies the Drinfel'd polynomial characterization of the usual definition of a KR module.

\begin{lemma}[{\cite[Lemma~3.4.1]{KKMMNN92}}]
\label{lemma:tensor_pairing}
Let $M_j$ and $N_j$, for $j = 1,2$, be $U_q'(\g)$-modules such that there exists a pairing $(\ ,\ )_j \colon M_j \times N_j \to \QQ(q)$ satisfying~\eqref{eq:admissible_pairing}. Then there exists a pairing $(\ ,\ ) \colon (M_1 \otimes M_2) \times (N_1 \otimes N_2) \to \QQ(q)$ defined by
\[
(u_1 \otimes u_2, v_1 \otimes v_2) = (u_1, v_1)_1 (u_2, v_2)_2,
\]
for all $u_j \in M_j$ and $v_j \in N_j$ with $j = 1,2$, that satisfies~\eqref{eq:admissible_pairing}.
\end{lemma}

\begin{proposition}[{\cite[Prop.~3.4.3]{KKMMNN92}}]
\label{prop:fusion_construction}
Let $u \in W(\varpi_r)_{\varpi_r}$ be a vector such that $\Absval{u}^2 = 1$.
\begin{enumerate}[\rm (1)]
\item The pairing $(\ ,\ ) \colon W^{r,s} \times W^{r,s} \to \QQ(q)$ constructed using Lemma~\ref{lemma:tensor_pairing} and the prepolarization on $W^{r,1}$ (see~\cite{Kashiwara02}) is a nondegenerate prepolarization on $W^{r,s}$.
\item $\Absval{ R_s(u^{\otimes s}) }^2 = 1$.
\item $\bigl( (W^{r,s})_{K_{\ZZ}}, (W^{r,s})_{K_{\ZZ}} \bigr) \subseteq K_{\ZZ}$, where
\[
(W^{r,s})_{K_{\ZZ}} = R_s\bigl( (U_q'(\g)_{K_{\ZZ}} u)^{\otimes s} \bigr) \cap \bigl( (U_q'(\g)_{K_{\ZZ}} u)^{\otimes s} \bigr)
\]
is a $U_q'(\g)_{K_{\ZZ}}$-submodule of $W^{r,s}$.
\end{enumerate}
\end{proposition}

\section{Existence of KR crystals}
\label{sec:existence}

This section is devoted to proving our main result.

\begin{theorem}
\label{thm:existence}
Let $r$ be such that $W^{r,s}$ is multiplicity free as a $U_q(\g_0)$-module for all $s \in \ZZ_{>0}$. Then $W^{r,s}$ admits a crystal pseudobase. Moreover, the KR crystal $B^{r,s}$ exists.
\end{theorem}

\begin{table}
\[
\begin{array}{cccccc}
\toprule
\g &  E_6^{(1)} & E_7^{(1)} & E_8^{(1)} & F_4^{(1)} & E_6^{(2)}
\\\midrule
r & 3,5 & 2,6 & 7 & 4 & 4
\\\bottomrule
\end{array}
\]
\caption{The nodes $r$ such that we show $B^{r,s}$ exists.}
\label{table:nodes}
\end{table}

We prove Theorem~\ref{thm:existence} case-by-case.
When $r$ is adjacent to $0$ or in the orbit of $0$ under a Dynkin diagram automorphism, Theorem~\ref{thm:existence} was shown in~\cite{KKMMNN92}. Theorem~\ref{thm:existence} was shown in nonexceptional affine types~\cite{Okado07,OS08}. Thus, it remains to show Theorem~\ref{thm:existence} for the values given in Table~\ref{table:nodes}.

From Proposition~\ref{prop:fusion_construction} and Proposition~\ref{prop:existence}, it is sufficient to show for the $U_q(\g_0)$-module decomposition $W^{r,s} \iso \bigoplus_{k=1}^M V(\overline{\lambda}_k)$ (where $\overline{\lambda}_k \in P_0^+$), there exists $u_k \in (M_{K_{\ZZ}})_{\lambda_k}$ such that
\begin{enumerate}[(i)]
\item \label{cond:norm} $(u_k, u_{\ell}) \in \delta_{k\ell} + q\mcA$ and
\item \label{cond:e} $\Absval{e_i u_k}^2 \in q_i^{-2\inner{h_i}{\lambda_k}-2} q \mcA$.
\end{enumerate}
The $U_q(\g_0)$-module decomposition of $W^{r,s}$ is given in~\cite{Chari01}.

We require the following facts.
Since the decomposition is multiplicity free, we have $(u_k, u_{\ell}) = 0$ for all $k \neq \ell$ since $\wt(u_k) \neq \wt(u_{\ell})$.
Note that
\[
[m] \in q^{1-m} \mcA,
\hspace{50pt}
\qbinom{m}{k}{q} \in q^{-k(m-k)} \mcA.
\]
Let $M$ be a $U_q'(\g)$-module.
We will use this variant of Equation~\eqref{eq:admissible_pairing}:
\begin{subequations}
\label{eq:divided_powers_admissible}
\begin{align}
\label{eq:divided_powers_admissible_e} (e_i^{(k)} v, w) & = q_i^{k(k-\inner{h_i}{\mu})} (v, f_i^{(k)} w),
\\ \label{eq:divided_powers_admissible_f} (f_i^{(k)} v, w) & = q_i^{k(k+\inner{h_i}{\mu})} (v, e_i^{(k)} w),
\end{align}
\end{subequations}
for all $w \in M_{\mu}$.
We also require
\begin{equation}
\label{eq:fe_relation}
f_i^{(a)} e_i^{(b)} v = \sum_{k=0}^{\min(a,b)} \qbinom{a-b-\inner{h_i}{\mu}}{k}{q_i} e_i^{(b-k)} f_i^{(a-k)} v,
\end{equation}
for any $v \in M_{\mu}$, which follows from applying the defining relation on $[e_i, f_i]$.
By applying Equation~\eqref{eq:divided_powers_admissible}, Equation~\eqref{eq:fe_relation}, and the bilinearity of $(\ ,\ )$, we have for any $v \in M_{\mu}$:
\begin{align*}
\Absval{e_i v}^2 & = q_i^{1-\inner{h_i}{\mu}} (v, f_i e_i v)
\\ & = q_i^{1-\inner{h_i}{\mu}} (v, e_i f_i v + [-\inner{h_i}{\mu}]_{q_i} v)
\\ & = q_i^{1-\inner{h_i}{\mu}} \bigl( (v, e_i f_i v) + [-\inner{h_i}{\mu}]_{q_i} (v, v) \bigr)
\\ & = q_i^{1-\inner{h_i}{\mu}} \left( q_i^{-(1+\inner{h_i}{\mu})} \Absval{f_i v}^2 + [-\inner{h_i}{\mu}]_{q_i} \Absval{v}^2 \right)
\end{align*}
Thus, we have
\begin{equation}
\label{eq:ei_norm}
\Absval{e_i v}^2 = q_i^{-2\inner{h_i}{\mu}} \Absval{f_i v}^2 + q_i^{1-\inner{h_i}{\mu}} [-\inner{h_i}{\mu}]_{q_i} \Absval{v}^2.
\end{equation}


For the remainder of the proof, we let $u \in W^{r,s}_{s\varpi_r}$ be such that $\Absval{u}^2 = 1$.
We have
\begin{equation}
\label{eq:fi_module_gen}
\Absval{f_i u}^2 = q_i^{1+\delta_{ir}s} (u, e_i f_i u) = q_i^{1+\delta_{ir} s} (u, [\delta_{ir} s]_{q_i} u) = q_i^{1+\delta_{ir}s} [\delta_{ir}s]_{q_i}
\end{equation}
for all $i \in I_0$ by Equation~\eqref{eq:divided_powers_admissible_e}, the defining relation on $[e_i, f_i]$ (or Equation~\eqref{eq:fe_relation}), and $e_i u = 0$. So we have $\Absval{f_r u}^2 \in q_r^2 \mcA$ (note $f_i u = 0$ for all $i \neq r$).

\subsection{Type \texorpdfstring{$E_6^{(1)}$, $r = 3$}{E6(1), r=3}}

We claim the elements
\[
u_k := e_6^{(k)} e_5^{(k)} e_4^{(k)} e_2^{(k)} e_0^{(k)} u
\]
are the desired elements, where $0 \leq k \leq s$. We have
\[
\wt(u_k) = \lambda_k := (s-k) \Lambda_3 + k\Lambda_6 - (2s-k)\Lambda_0,
\]
and from~\cite{Chari01}, the classical decomposition is $W^{3,s} \iso \bigoplus_{k=0}^s V\bigl( (s-k) \clfw_3 + k \clfw_6 \bigr)$. Thus, we need to show $u_k$ satisfies~(\ref{cond:norm}) and~(\ref{cond:e}).

We first show~(\ref{cond:norm}). We have
\[
\Absval{u_k}^2 = q_6^{k(k-k)} (e_5^{(k)} e_4^{(k)} e_2^{(k)} e_0^{(k)} u, f_6^{(k)} u_k)
\]
from Equation~\eqref{eq:divided_powers_admissible_e}. Next, we have
\begin{equation}
\label{eq:commute_f}
\begin{split}
f_6^{(k)} u_k & = f_6^{(k)} e_6^{(k)} e_5^{(k)} e_4^{(k)} e_2^{(k)} e_0^{(k)} u
\\ & = \sum_{m=0}^k \qbinom{k}{m}{q_6} e_6^{(k-m)} f_6^{(k-m)} e_5^{(k)} e_4^{(k)} e_2^{(k)} e_0^{(k)} u
\\ & = e_5^{(k)} e_4^{(k)} e_2^{(k)} e_0^{(k)} u,
\end{split}
\end{equation}
where the second equality comes from Equation~\eqref{eq:fe_relation} and the third equality follows from the fact $e_i f_j = f_j e_i$ for all $i \neq j$ and $f_6 w = 0$ (so only the $m = k$ term is nonzero).
By computations similar to Equation~\eqref{eq:commute_f}, we have
\[
\Absval{u_k}^2 = (e_5^{(k)} e_4^{(k)} e_2^{(k)} e_0^{(k)} u, e_5^{(k)} e_4^{(k)} e_2^{(k)} e_0^{(k)} u) = \Absval{e_0^{(k)} u}^2.
\]
Moreover, similar to Equation~\eqref{eq:commute_f}, we have
\begin{align*}
\Absval{e_0^{(k)} u}^2 & = (e_0^{(k)} u, e_0^{(k)} u) = q_0^{k(k+2s-2k)} (u, f_0^{(k)} e_0^{(k)} u)
\\ & = q_0^{k(2s-k)} \sum_{m=0}^k \qbinom{2s}{m}{q_0} (u, e_0^{(k-m)} f_0^{(k-m)} u) = q_0^{k(2s-k)} \qbinom{2s}{k}{q_0} (u, u)
\end{align*}
since $f_0 u = 0$. Hence, we have
\begin{equation}
\label{eq:norm_uk}
\Absval{u_k}^2 = q_0^{k(2s-k)} \qbinom{2s}{k}{q_0} \in 1 + q\mcA.
\end{equation}

Next, we show~(\ref{cond:e}). Fix some $i \in I_0$. From Equation~\eqref{eq:ei_norm}, it remains to compute $\Absval{f_i u_k}^2$.
We compute $\Absval{f_i u_k}^2$ depending on the value of $i$. We note that the case of $k = 0$ is done by Equation~\eqref{eq:fi_module_gen}. Therefore, we assume $k \geq 1$. For $i = 6$, we have
\begin{equation}
\label{eq:f_first}
f_6 u_k = \qbinom{1-k+k}{1}{q_6} e_6^{(k-1)} e_5^{(k)} e_4^{(k)} e_2^{(k)} e_0^{(k)} u + e_6^{(k)} f_6 e_5^{(k)} e_4^{(k)} e_2^{(k)} e_0^{(k)} u = e_6^{(k-1)} e_5^{(k)} e_4^{(k)} e_2^{(k)} e_0^{(k)} u
\end{equation}
by Equation~\eqref{eq:fe_relation} and the fact $f_6 u = 0$.
Hence, similar to the computation for $\Absval{u_k}^2$, we have
\begin{align*}
\Absval{f_6 u_k}^2 & =  q_6^{k-1} \Absval{e_6^{(k-1)} e_5^{(k)} e_4^{(k)} e_2^{(k)} e_0^{(k)} u}^2
\\ & =  q_6^{k-1} \qbinom{k}{k-1}{q_6} \Absval{e_5^{(k)} e_4^{(k)} e_2^{(k)} e_0^{(k)} u}^2
\\ & = q_6^{k-1} \qbinom{k}{k-1}{q_6} q_0^{k(2s-k)} \qbinom{2s}{k}{q_0}.
\end{align*}
For $i = 1$, we have $f_1 u_k  = e_6^{(k)} e_5^{(k)} e_4^{(k)} e_2^{(k)} e_0^{(k)} f_1 u = 0$, and so $\Absval{f_1 u_k}^2 = 0$.
For $i = 5,4,2$, we have $f_i u_k = 0$ by applying Equation~\eqref{eq:fe_relation} and the Serre relations (\textit{e.g.}, a straightforward calculation shows $e_4^{(k)} e_2^{(k-1)} e_0^{(k)} u = 0$ by repeatedly applying the Serre relations).
Finally, we have $f_3 u_k = e_6^{(k)} e_5^{(k)} e_4^{(k)} e_2^{(k)} e_0^{(k)} f_3 u$. Therefore, we have $\Absval{f_3 u_k}^2 = \Absval{e_4^{(k)} e_2^{(k)} e_0^{(k)} f_3 u}^2$ similar to Equation~\eqref{eq:commute_f}. However, for removing $e_4^{(k)}$, we obtain
\[
(e_4^{(k)} e_2^{(k)} e_0^{(k)} f_3 u, e_4^{(k)} e_2^{(k)} e_0^{(k)} f_3 u) = q_4^{k(k-(k+1))} (e_2^{(k)} e_0^{(k)} f_3 u, f_4^{(k)} e_4^{(k)} e_2^{(k)} e_0^{(k)} f_3 u)
\]
by Equation~\eqref{eq:divided_powers_admissible_e}. Furthermore, we have
\begin{align*}
f_4^{(k)} e_4^{(k)} e_2^{(k)} e_0^{(k)} f_3 u & = \sum_{m=0}^k \qbinom{k-1}{m}{q_4} e_4^{(k-m)} f_4^{(k-m)} e_2^{(k)} e_0^{(k)} f_3 u
\\ & = \qbinom{k-1}{k-1}{q_4} e_4 e_2^{(k)} e_0^{(k)} f_4 f_3 u + \qbinom{k-1}{k}{q_4} e_2^{(k)} e_0^{(k)} f_3 u
\\ & = e_4 e_2^{(k)} e_0^{(k)} f_4 f_3 u,
\end{align*}
where we note that $\qbinom{k-1}{k}{q_4} = 0$ (recall that we assumed $k \geq 1$).
Thus, by applying Equation~\eqref{eq:divided_powers_admissible_e}, we obtain
\begin{equation}
\label{eq:moving_e4}
\Absval{e_4^{(k)} e_2^{(k)} e_0^{(k)} f_3 u}^2 = q_4^{-k} (e_2^{(k)} e_0^{(k)} f_3 u, e_4 e_2^{(k)} e_0^{(k)} f_4 f_3 u) = q_4^{-k} q_4^k \Absval{e_2^{(k)} e_0^{(k)} f_4 f_3 u}^2.
\end{equation}
Next, we have that
\[
\Absval{e_2^{(k)} e_0^{(k)} f_4 f_3 u}^2 = \Absval{ e_0^{(k)} f_2 f_4 f_3 u}^2
\]
from a similar computation to Equation~\eqref{eq:moving_e4}.
Continuing using Equation~\eqref{eq:divided_powers_admissible_e}, we have
\[
\Absval{ e_0^{(k)} f_2 f_4 f_3 u}^2 = q_0^{k(2s-1-k)} (f_2 f_4 f_3 u, f_0^{(k)} e_0^{(k)} f_2 f_4 f_3 u).
\]
We note that $f_0 f_2 f_4 f_3 w = 0$ for any $w \in W^{3,1}_{\varpi_3}$ from weight considerations and the classical decomposition. So $f_0 f_2 f_4 f_3 (w_1 \otimes \cdots \otimes w_s) = 0$ for any $w_1, \dotsc, w_s \in W^{3,1}_{\varpi_3}$ from applying the coproduct $\Delta(f_i) = f_i \otimes 1 + K_i \otimes f_i$. Thus, we have $f_0 f_2 f_4 f_3 u = 0$ from the construction of $u$ and $W^{3,s}$.
Therefore, we compute
\[
f_0^{(k)} e_0^{(k)} f_2 f_4 f_3 u = \sum_{m=0}^k \qbinom{2s-1}{m}{q_0} e_0^{(k-m)} f_0^{(k-m)} f_2 f_4 f_3 u =\qbinom{2s-1}{k}{q_0} f_2 f_4 f_3 u
\]
similar to Equation~\eqref{eq:commute_f} and using the Serre relations.
Thus, we have
\[
\Absval{e_0^{(k)} f_2 f_4 f_3 u}^2 = q_0^{k(2s-1-k)} \qbinom{2s-1}{k}{q_0} \Absval{f_2 f_4 f_3 u}^2.
\]
Next, we see
\begin{align*}
\Absval{f_2 f_4 f_3 u}^2 & = q_2^{1-1} (f_4 f_3 u, e_2 f_2 f_4 f_3 u) = (f_4 f_3 u, [1]_{q_2} f_4 f_3 u)
\\ & = q_4^{1-1} (f_3 u, e_4 f_4 f_3 u) = (f_3 u, [1]_{q_4} f_3 u)  = \Absval{f_3 u}^2
\end{align*}
by a similar computation to Equation~\eqref{eq:fi_module_gen}.
Hence, we have
\begin{equation}
\label{eq:f_r_hw}
\Absval{f_3 u_k}^2 = q_0^{k(2s-1-k)} \qbinom{2s-1}{k}{q_0} \Absval{f_3 u}^2 = q_0^{k(2s-1-k)} \qbinom{2s-1}{k}{q_0} q_3^{1+s} [s]_{q_3} \in q_3^2 \mcA,
\end{equation}
where the last equality is by Equation~\eqref{eq:fi_module_gen}.
To complete the proof of~(\ref{cond:e}), we can see that
\[
q_i^{-2\inner{h_i}{\lambda_k}} \Absval{f_i u_k}^2 \in q_i^{-2\inner{h_i}{\lambda_k}} \mcA,
\qquad
q_i^{1-\inner{h_i}{\lambda_k}} [-\inner{h_i}{\lambda_k}]_{q_i} q_0^{k(2s-1-k)} \qbinom{2s-1}{k}{q_0} \in q_i^2 \mcA,
\]
noting $\inner{h_i}{\lambda_k} \geq 0$.

\subsection{Type \texorpdfstring{$E_6^{(1)}$, $r = 5$}{E6(1), r=5}}

The following are the desired elements in $W^{5,s}$:
\[
u_k := e_1^{(k)} e_3^{(k)} e_4^{(k)} e_2^{(k)} e_0^{(k)} u_0 \in W^{5,s}_{(s-k) \varpi_5 + k \varpi_1},
\]
where $0 \leq k \leq s$.
The proof is the same as $r = 3$ after applying the order $2$ diagram automorphism that fixes $0$.

\subsection{Type \texorpdfstring{$E_7^{(1)}$, $r = 2$}{E7(1), r=2}}

The following are the desired elements in $W^{2,s}$:
\[
u_k := e_7^{(k)} e_6^{(k)} e_5^{(k)} e_4^{(k)} e_3^{(k)} e_1^{(k)} e_0^{(k)} u_0 \in W^{2,s}_{(s-k) \varpi_2 + k \varpi_7},
\]
where $0 \leq k \leq s$.
The proof is similar to the $W^{3,s}$ in type $E_6^{(1)}$, where we compute
\begin{align*}
\Absval{u_k}^2 & = q_0^{k(2s-k)} \qbinom{2s}{k}{q_0},
\\ \Absval{f_7 u_k}^2 & =  q_7^{k-1} \qbinom{k}{k-1}{q_7} \Absval{u_k}^2,
\\ \Absval{f_i u_k}^2 & = 0 \hspace{50pt} (i= 6,5,4,3,1),
\\ \Absval{f_2 u_k}^2 & = q_0^{k(2s-1-k)} \qbinom{2s-1}{k}{q_0} \Absval{f_2 u}^2. 
\end{align*}

\subsection{Type \texorpdfstring{$E_6^{(2)}$, $r = 4$}{E6(2), r=4}}

We claim
\[
u_{k',k} := e_0^{(k')} e_1^{(k)} e_2^{(k)} e_3^{(k)} e_2^{(k)} e_1^{(k)} e_0^{(k)} u
\]
are the desired elements, where $0 \leq k' \leq k \leq s$.
We note that
\[
\wt(u_{k',k}) = \lambda_{k',k} := (s-k)\Lambda_4 + (k-k')\Lambda_1 - (2s-2k')\Lambda_0.
\]
To obtain the parameterization of the classical decomposition
\[
W^{4,s} \iso \bigoplus_{\substack{t_1,t_2 \geq 0 \\ t_1 + t_2 \leq s}} V(t_1 \clfw_4 + t_2 \clfw_1)
\]
given in~\cite[Prop.~9.31]{Scrimshaw17}, we set $t_1 = s - k$ and $t_2 = k - k'$ (which is forced by weight considerations).
Note that $t_1 \geq 0$ if and only if $k \leq s$; $t_2 \geq 0$ if and only if $k' \leq k$; and $t_1 + t_2 \leq s$ if and only if $0 \leq k'$ (as $t_1 + t_2 = s - k'$).
Hence, we have the same classical decomposition.

To show~(\ref{cond:norm}), we have
\[
\Absval{u_{0,k}}^2 = q_1^{k(k-k)} ( e_2^{(k)} e_3^{(k)} e_2^{(k)} e_1^{(k)} e_0^{(k)} u, f_1^{(k)} u_{0,k} ).
\]
Next, we compute
\begin{align*}
f_1^{(k)} u_{0,k} & = f_1^{(k)} e_1^{(k)} e_2^{(k)} e_3^{(k)} e_2^{(k)} e_1^{(k)} e_0^{(k)} u
\\ & = \sum_{m=0}^{k} \qbinom{k}{m}{q_1} e_1^{(k-m)} f_1^{(k-m)} e_2^{(k)} e_3^{(k)} e_2^{(k)} e_1^{(k)} e_0^{(k)} u
\\ & = \sum_{m=0}^{k} \qbinom{k}{m}{q_1} e_1^{(k-m)} e_2^{(k)} e_3^{(k)} e_2^{(k)} \sum_{p=0}^{k-m} \qbinom{k-m}{p}{q_1} e_1^{(k-p)} f_1^{(k-m-p)} e_0^{(k)} u
\\ & = \sum_{m=0}^{k} \qbinom{k}{m}{q_1} \qbinom{k-m}{k-m}{q_1} e_1^{(k-m)} e_2^{(k)} e_3^{(k)} e_2^{(k)} e_1^{(m)} e_0^{(k)} u
\\ & = e_2^{(k)} e_3^{(k)} e_2^{(k)} e_1^{(k)} e_0^{(k)} u,
\end{align*}
where the last equality follows from the fact $e_2^{(k)} e_1^{(m)} e_0^{(k)} u = 0$ for all $k > m$ by the Serre relations and $e_2 u = 0$.
Hence, we have
\[
\Absval{u_{0,k}}^2 = \Absval{e_2^{(k)} e_3^{(k)} e_2^{(k)} e_1^{(k)} e_0^{(k)} u}^2 = q^{k(k-k)} (e_3^{(k)} e_2^{(k)} e_1^{(k)} e_0^{(k)} u, f_2^{(k)} e_2^{(k)} e_3^{(k)} e_2^{(k)} e_1^{(k)} e_0^{(k)} u).
\]
Now, similar to the previous computation, we obtain
\begin{align*}
f_2^{(k)} e_2^{(k)} e_3^{(k)} e_2^{(k)} e_1^{(k)} e_0^{(k)} u & = \sum_{m=0}^{k} \qbinom{k}{m}{q_2} e_2^{(k-m)} f_2^{(k-m)} e_3^{(k)} e_2^{(k)} e_1^{(k)} e_0^{(k)} u
\\ & = \sum_{m=0}^{k} \qbinom{k}{m}{q_2} e_2^{(k-m)} e_3^{(k)} \sum_{p=0}^{k-m} \qbinom{k-m}{p}{q_2} e_2^{(k-p)} f_2^{(k-m-p)} e_1^{(k)} e_0^{(k)} u
\\ & = \sum_{m=0}^{k} \qbinom{k}{m}{q_2} \qbinom{k-m}{k-m}{q_2} e_2^{(k-m)} e_3^{(k)} e_2^{(m)} e_1^{(k)} e_0^{(k)} u
\\ & = e_3^{(k)} e_2^{(k)} e_1^{(k)} e_0^{(k)} u
\end{align*}
since $e_3^{(k)} e_2^{(m)} e_1^{(k)} e_0^{(k)} u = 0$ for all $k > m$ by the Serre relations (recall that $A_{32} = -1$) and $e_3 u = 0$. Hence, we have
\[
\Absval{u_{0,k}}^2 = \Absval{e_3^{(k)} e_2^{(k)} e_1^{(k)} e_0^{(k)} u}^2 = q_0^{k(2s-k)} \qbinom{2s}{k}{q_0} \in 1 + q \mcA,
\]
where the last equality is shown similar to Equation~\eqref{eq:norm_uk}.

Next, we consider
\[
\Absval{u_{k',k}}^2 = q_0^{k'(k'+2s-2k')} (u_{0,k}, f_0^{(k')} u_{k', k}).
\]
We compute
\begin{equation}
\label{eq:f0_uk1k2}
f_0^{(k')} u_{k', k} = f_0^{(k')} e_0^{(k')} u_{0,k} = \sum_{m=0}^{k'} \qbinom{2s}{m}{q_0} e_0^{(k'-m)} f_0^{(k'-m)} u_{0,k},
\end{equation}
and
\[
f_0^{(k'-m)} e_0^{(k)} u = \sum_{p=0}^{k'-m} \qbinom{k'-m-k+2s}{p}{q_0} e_0^{(k-p)} f_0^{(k'-m-p)} u = \qbinom{k'-m-k+2s}{k'-m}{q_0} e_0^{(k-k'+m)} u
\]
as $k' - m \leq k$ (since $k' \leq k$ and $m \geq 0$) and $f_0 u = 0$. Next, we have $e_1^{(k)} e_0^{(m)} u = 0$ for all $k > m$ by the Serre relations and $e_1 u = 0$, and so the only term that is nonzero in Equation~\eqref{eq:f0_uk1k2} is when $m = k'$. Therefore, we have
\[
\Absval{u_{k',k}}^2 = q_0^{k'(2s-k')} \qbinom{2s}{k'}{q_0} \Absval{u_{0,k}}^2 = q_0^{k'(2s-k')} \qbinom{2s}{k'}{q_0} q_0^{k(2s-k)} \qbinom{2s}{k}{q_0} \in 1 + q \mcA.
\]

To show~(\ref{cond:e}), it remains to compute $\Absval{f_i u_{k',k}}^2$ by Equation~\eqref{eq:ei_norm}, and by Equation~\eqref{eq:fi_module_gen}, we can assume $k \geq 1$. For $i \in I_0$, we have $f_i u_{k',k} = e_0^{(k')} f_i u_{0,k}$, and by the above, we have
\[
\Absval{f_i u_{k',k}}^2 = q_0^{k'(2s-\delta_{i1}-k')} \qbinom{2s-\delta_{i1}}{k'}{q_0} \Absval{f_i u_{0,k}}^2.
\]
Next, similar to the computation in Equation~\eqref{eq:f_first}, we have
\begin{align*}
f_1 u_{0,k} & = e_1^{(k-1)} e_2^{(k)} e_3^{(k)} e_2^{(k)} e_1^{(k)} e_0^{(k)} u + e_1^{(k)} e_2^{(k)} e_3^{(k)} e_2^{(k)} f_1 e_1^{(k)} e_0^{(k)} u
\\ & = e_1^{(k-1)} e_2^{(k)} e_3^{(k)} e_2^{(k)} e_1^{(k)} e_0^{(k)} u + e_1^{(k)} e_2^{(k)} e_3^{(k)} e_2^{(k)} e_1^{(k-1)} e_0^{(k)} u
\\ & = e_1^{(k-1)} e_2^{(k)} e_3^{(k)} e_2^{(k)} e_1^{(k)} e_0^{(k)} u,
\end{align*}
where the last equality is using $e_2^{(k)} e_1^{(m)} e_0^{(k)} u = 0$ for all $k > m$.
Therefore, we have
\[
\Absval{f_1 u_{0,k}}^2 = q_1^{k-1} \qbinom{k}{k-1}{q_1} q_0^{k(2s-k)} \qbinom{2s}{k}{q_0}
\]
by a computation similar to Equation~\eqref{eq:norm_uk}. Similar to Equation~\eqref{eq:f_r_hw}, we have
\[
\Absval{f_4 u_{0,k}}^2 = q_0^{k(2s-1-k)} \qbinom{2s-1}{k}{q_0} \Absval{f_4 u}^2.
\]
We also have $f_2 u_{0,k} = f_3 u_{0,k} = 0$ by applying the Serre relations.
Thus, we see that~(\ref{cond:e}) holds.

\subsection{Type \texorpdfstring{$E_7^{(1)}$, $r = 6$}{E7(1), r=6}}

The following are the desired elements in $W^{6,s}$:
\[
u_{k',k} := e_0^{(k')} e_1^{(k)} e_3^{(k)} e_4^{(k)} e_5^{(k)} e_2^{(k)} e_4^{(k)} e_3^{(k)} e_1^{(k)} e_0^{(k)} u \in W^{6,s}_{(s-t_1-t_2) \varpi_6 + t_2 \varpi_1},
\]
where $0 \leq k' \leq k \leq s$.
Then $\wt(u_{k',k}) = (s-k)\Lambda_6 + (k-k')\Lambda_1 - (2s-2k')\Lambda_0$.
Showing the classical decomposition is the same as in~\cite{Chari01} is similar to the $r=4$ case for type $E_6^{(2)}$.
Moreover, it is similar to show that
\begin{align*}
\Absval{u_{k',k}}^2 & = q_0^{k'(2s-k')} \qbinom{2s}{k'}{q_0} q_0^{k(2s-k)} \qbinom{2s}{k}{q_0}, \allowdisplaybreaks
\\ \Absval{f_i u_{k',k}}^2 & = q_0^{k'(2s-\delta_{i1}-k')} \qbinom{2s-\delta_{i1}}{k'}{q_0} \Absval{f_i u_{0,k}}^2 \hspace{50pt} (i \in I_0), \allowdisplaybreaks
\\ \Absval{f_1 u_{0,k}}^2 & = q_1^{k-1} \qbinom{k}{k-1}{q_1} \Absval{u_{0,k}}^2,
\\ \Absval{f_i u_{0,k}}^2 & = 0 \hspace{50pt} (i=2,3,4,5,7),
\\ \Absval{f_6 u_{0,k}}^2 & = q_0^{k(2s-1-k)} \qbinom{2s-1}{k}{q_0} \Absval{f_6 u}^2.
\end{align*}

\subsection{Type \texorpdfstring{$E_8^{(1)}$, $r = 1$}{E8(1), r=1}}

The following are the desired elements in $W^{1,s}$:
\[
u_{k',k} := e_0^{(k')} e_8^{(k)} e_7^{(k)} e_6^{(k)} e_5^{(k)} e_4^{(k)} e_3^{(k)} e_2^{(k)} e_4^{(k)} e_5^{(k)} e_6^{(k)} e_7^{(k)} e_8^{(k)} e_0^{(k)} u \in W^{1,s}_{(s-t_1-t_2) \varpi_1 + t_2 \varpi_8},
\]
where $0 \leq k' \leq k \leq s$.
Then $\wt(u_{k',k}) = (s-k)\Lambda_1 + (k-k')\Lambda_8 - (2s-2k')\Lambda_0$.
Showing the classical decomposition is the same as in~\cite{Chari01} is similar to the $r=4$ case for type $E_6^{(2)}$.
Moreover, it is similar to show that
\begin{align*}
\Absval{u_{k',k}}^2 & = q_0^{k'(2s-k')} \qbinom{2s}{k'}{q_0} q_0^{k(2s-k)} \qbinom{2s}{k}{q_0}, \allowdisplaybreaks
\\ \Absval{f_i u_{k',k}}^2 & = q_0^{k'(2s-\delta_{i8}-k')} \qbinom{2s-\delta_{i8}}{k'}{q_0} \Absval{f_i u_{0,k}}^2 \hspace{50pt} (i \in I_0), \allowdisplaybreaks
\\ \Absval{f_8 u_{0,k}}^2 & = q_8^{k-1} \qbinom{k}{k-1}{q_8} \Absval{u_{0,k}}^2,
\\ \Absval{f_i u_{0,k}}^2 & = 0 \hspace{50pt} (i=2,3,4,5,6,7),
\\ \Absval{f_1 u_{0,k}}^2 & = q_0^{k(2s-1-k)} \qbinom{2s-1}{k}{q_0} \Absval{f_1 u}^2.
\end{align*}

\subsection{Type \texorpdfstring{$F_4^{(1)}$, $r = 4$}{F4(1), r=4}}

The following are the desired elements in $W^{4,s}$:
\[
u_{k',k} := e_0^{(k')} e_1^{(k)} e_2^{(k)} e_3^{(2k)} e_2^{(k)} e_1^{(k)} e_0^{(k)} u \in W^{4,s}_{(s-2k) \varpi_4 + (k - k') \varpi_1},
\]
where $0 \leq k' \leq k \leq s/2$.
Then $\wt(u_{k',k}) = (s-2k)\Lambda_4 + (k-k')\Lambda_1 - (s-2k')\Lambda_0$.
To obtain the parameterization of the classical decomposition
\[
W^{4,s} \iso \bigoplus_{t_2=0}^{s/2} \bigoplus_{t_1=0}^{t_2} V\bigl( (s-2t_2) \clfw_4 + t_1 \clfw_1 \bigr)
\]
given in~\cite{Chari01}, we take $t_1 = k - k'$ and $t_2 = k$. Indeed, we have $t_2 \leq s/2$ if and only if $k \leq s/2$; $t_1 \geq 0$ if and only if $k \leq k'$; and $t_1 \leq t_2$ if and only if $0 \leq k'$.

Moreover, it is similar to the $r=4$ case for type $E_6^{(2)}$ to show that
\begin{align*}
\Absval{u_{k',k}}^2 & = q_0^{k'(2s-k')} \qbinom{2s}{k'}{q_0} q_0^{k(2s-k)} \qbinom{2s}{k}{q_0}, \allowdisplaybreaks
\\ \Absval{f_i u_{k',k}}^2 & = q_0^{k'(2s-\delta_{i1}-k')} \qbinom{2s-\delta_{i1}}{k'}{q_0} \Absval{f_i u_{0,k}}^2 \hspace{50pt} (i \in I_0), \allowdisplaybreaks
\\ \Absval{f_1 u_{0,k}}^2 & = q_1^{k-1} \qbinom{k}{k-1}{q_1} \Absval{u_{0,k}}^2,
\\ \Absval{f_2 u_{0,k}}^2 & = \Absval{f_3 u_{0,k}}^2 = 0,
\\ \Absval{f_4 u_{0,k}}^2 & = q_0^{k(2s-1-k)} \qbinom{2s-1}{k}{q_0} \Absval{f_4 u}^2.
\end{align*}

\bibliographystyle{alpha}
\bibliography{existence}{}
\end{document}